\documentclass[11pt]{article}

\usepackage{amssymb,amsmath,color}

\newtheorem{Theorem}{Theorem}[section]

\newtheorem{Lemma}{Lemma}[section]
\newtheorem{Corollary}{Corollary}[section]

\numberwithin{equation}{section}

\newcommand{\non}{\nonumber}
\newcommand{\cR}{\mathbb R}
\newcommand{\cN}{{\rm I\!N}}

\newcommand{\f}{\frac}

\newcommand{\eq}[1]{\mbox{\rm {(\ref{#1})}}}

\title{\Large\bf Spherically symmetric solutions to a model for phase\\
 transitions driven by configurational forces
}

\author{\small
\sc Yaobin Ou$^{1}$\thanks{E-mail: ou@bcamath.org} and
\small \sc Peicheng Zhu$^{1,2}$\thanks{E-mail: zhu@bcamath.org}\\
\small $^1$ Basque Center for Applied Mathematics (BCAM)\\
\small Building 500, Bizkaia Technology Park \\
\small E-48160 Derio, \ \ Spain\\
\small $^2$ IKERBASQUE, Basque Foundation for Science\\
\small E-48011 Bilbao,\ \ Spain}

\date{ }

\begin{document}

\maketitle

\begin{abstract}

We prove the global in time existence of spherically symmetric solutions  to an initial-boundary value
problem for a system of partial differential equations, which consists of the
equations of linear elasticity and a nonlinear, non-uniformly parabolic
equation of second order. The problem models the behavior in time of materials
 in which martensitic phase transitions, driven by configurational forces, take place, and can be considered to be a regularization
of the corresponding sharp interface model. By assuming that the solutions are spherically symmetric,
 we reduce the original multidimensional problem to the one in one space dimension, then prove the existence of spherically symmetric
 solutions.
   Our proof is valid due to the essential feature that the reduced problem is one space dimensional.

\end{abstract}
\section{Introduction}

Many inhomogeneous systems can be characterized by domains of different phases
separated by a distinct interface. When driven out of
equilibrium, their dynamics result in the evolution of those interfaces, and the systems
might develop into structures (compositional and structural inhomogeneities) with characteristic length scales at the
nano-, micro- or meso-scale. To a large extent, the material properties of such
systems are determined by those structures of small-scale. Thus it is  important to
understand precisely the mechanisms that drive the evolution of those structures.
 Materials microstructures may consist of spatially distributed phases of different
compositions and/or crystal structures, grains of different orientations, domains of
different structural variants, domains of different electrical or magnetic polarizations,
and structural defects. These structural features usually have an intermediate
mesoscopic length scale in the range of nanometers to microns. The size, shape,
and spatial arrangement of the local structural features in a microstructure play a
critical role in determining the physical properties of a material. Because of the complex
and nonlinear nature of microstructure evolution, numerical approaches are often
employed.  For more details, see e.g. \cite{Chen02,Emmerich03,Moelans}.

 In this article we are interested in a model for the evolution, driven
 by configurational forces, of microstructures in  elastically deformable solids. There are two main types of modeling for the evolution of  microstructures.
 In the conventional approach, the regions separating the domains are treated as mathematically
sharp interfaces. The local interfacial velocity is then determined as part
of the boundary conditions, or is calculated from the driving force for interface motion
and the interfacial mobility. This approach requires the explicit tracking of the interface
positions. Such an interface-tracking approach can be successful in one-dimensional
systems, however it will be impractical for complicated three-dimensional microstructures.
Therefore, during the past decades, another approach
 has been invented, namely, the phase-field approach in which the interface is  not of zero
  thickness, instead an interfacial region with thickness of certain order of a small regularization
   parameter. Though it is still a young discipline in condensed matter physics,  this approach has emerged to be
 one of the most powerful methods for modeling the evolution of microstructures.
 It can be traced back  the theory of  diffuse-interface description, which is developed, independently, more
than a century ago by van der Waals \cite{VanderWaals} and some half century ago
by Cahn and Hilliard \cite{Cahn58}.

The two well-known models for temporal evolution of microstructures are
 the  Cahn-Hilliard/Allen-Cahn  equations corresponding, respectively, to the case that the order parameter is conserved and not
 conserved. These phase field models describe microstructure phenomena at the
mesoscale (see e.g. \cite{Moelans}), and one suitable limit of it may be the corresponding sharp-  or thin-interface descriptions.
 In this article we study a model for the behavior in time of materials with
diffusionless phase transitions. The model has diffusive interfaces and
consists of the partial differential equations of linear elasticity coupled to
a quasilinear, non-uniformly parabolic equation of second order that differs from the Allen-Cahn equation
 (the Cahn-Hilliard equation in the case that the order parameter is conserved) by a
gradient term. It is derived in \cite{Alber00,Alber04} from a sharp interface model for diffusionless phase
transitions and can be considered to be a regularization of that model.
 To verify the validity of the new model,  mathematical analysis has been carried
 out for the existence/regularity of weak solutions to initial boundary value problems in
 one space dimension, \cite{Alber06,Alber07,Alber10a,Zhu1,Zhu2}, the motion of interfaces \cite{Alber09},
 and the existence of traveling waves \cite{Kawashima10}. In the present article, the existence
 of spherically symmetric solutions to an initial-boundary value problem will be studied.
 We first formulate this initial-boundary value problem in the
three-dimensional case, then reduce it,  by assuming that the solution is spherically symmetric,
to the one-dimensional case. The existence of weak solutions to this one dimensional
problem is proved.

\vskip0.2cm
Let $\Omega\subset\cR^3 $ be an open set. It represents
the material points of a solid body. The different phases are
characterized by the order parameter $S(t,x)\in \cR$. A value of
$S(t,x)$ close to zero indicates that the material is in the matrix
phase at the point $x\in\Omega$ at time $t$, a value close to one
indicates that the material is in the second phase. The other unknowns
are the displacement $u(t,x)\in \cR^3 $ of the material point $x$ at
time $t$ and the Cauchy stress tensor $T(t,x)\in {\cal S}^3 $, where
${\cal S}^3 $ denotes the set of symmetric $3\times 3$-matrices. The
unknowns must satisfy the quasi-static equations
\begin{eqnarray}
 -{\rm div}_x\,T(t,x)&=&b(t,x),
 \label{eq1}  \\
 T(t,x) &=& D \big(\varepsilon(\nabla_x\, u(t,x)) - \bar\varepsilon S(t,x) \big),
 \label{eq1a} \\
 S_t(t,x) &=& -c \Big( \psi_S(\varepsilon(\nabla_x
            u(t,x)),S(t,x)) - \nu \Delta_x S(t,x) \Big)|\nabla_x S(t,x)| \quad
 \label{eq2}
\end{eqnarray}
for $(t,x)  \in (0,\infty)\times\Omega$. The boundary and initial
conditions are
\begin{eqnarray}
 u(t,x)=\gamma(t,x),&& S(t,x)=0, \hspace{46pt} (t,x)\in
 [0,\infty)\times\partial\Omega,
 \label{eq3a} \\
 && S(0,x) = S_0(x),\qquad x\in \Omega.
 \label{eq3}
\end{eqnarray}
Here $\nabla_xu$ denotes the $3\times 3$-matrix of first order
derivatives of $u$, the deformation gradient, $(\nabla_xu)^T$ denotes the
transposed matrix and
$$
 \varepsilon(\nabla_xu) = \frac12\left(\nabla_x u + (\nabla_x u ) ^T \right)
$$
is the strain tensor. $\bar\varepsilon\in {\cal S}^3 $ is a given matrix,
the misfit strain, and $D:{\cal S}^3\to {\cal S}^3$ is the elasticity
tensor, a linear, symmetric, positive definite mapping. In the free energy
\begin{equation}\label{freengy}
 \psi(\varepsilon,S)=\f12 \big(D(\varepsilon-\bar\varepsilon S)\big)
 \cdot (\varepsilon-\bar\varepsilon S)+\hat\psi(S),
\end{equation}
we assume that $\hat\psi\in C^2(\cR,[0,\infty))$,  choose $\hat\psi$
 as a double well potential with
minima at $S=0$ and $S=1$. $\psi_S$ is the partial derivative.
The scalar product of two matrices $A$ and $B$
is denoted by $A \cdot B =
\sum a_{ij} b_{ij}$.  $c>0$ is a
constant and $\nu$ is a small positive constant. Given are the volume force
$b:[0,\infty)\times\Omega\to \cR^3$ and the data $\gamma:[0,\infty)\times
\partial\Omega\to \cR^3$, $S_0: \Omega\to \cR$.

This completes the formulation of the initial-boundary value problem.
Equations (\ref{eq1}) and (\ref{eq1a}) differ from the system of linear
elasticity only by the term $\bar\varepsilon S$. The evolution equation
(\ref{eq2}) for the order parameter $S$ is non-uniformly parabolic because of
the term $\nu\Delta S |\nabla_xS|$. Since this initial-boundary value problem
is derived from a sharp interface model, to verify that it is indeed a
diffusive interface model regularizing the sharp interface model, it must be
shown that the equations (\ref{eq1}) -- (\ref{eq3}) with positive $\nu$ have
solutions which exist globally in time, and that these solutions tend to
solutions of the sharp interface model for $\nu\to 0$. This would also be a
method to prove existence of solutions to the original sharp interface model.

We only contribute to the first part of this program in this work and show
that there exist some special solutions to  the initial-boundary value problem that is
essentially in one space dimension. Up to now we still can't solve the following problem:
 either solutions in three space dimensions exist or  these
solutions converge to a solution of the sharp interface model for $\nu
\rightarrow 0$. We shall see later that the existence result
of spherically symmetric solutions is of interest because
 the problem has a stronger nonlinear term (compared with
 the problem by assuming all unknowns depend on
 only one component of space variable $x$ which is studied in \cite{Alber06}),
 despite it is essentially one space dimensional.

Related to our investigations is the model for diffusion dominated phase
transformations obtained by coupling the elasticity equations (\ref{eq1}),
(\ref{eq1a}) with the Allen-Cahn/Cahn-Hilliard equations. They  have recently been studied in \cite{BCDGSS02, CMP00,
Garcke03}.

\bigskip
\noindent {\bf Statement of the main result.} Since we shall look for solutions, which are spherically symmetric,
 to problem \eq{eq1} -- \eq{eq3},  the problem can be reduced to the one which is one space dimensional. To this end
we now assume that the body force boundary and initial data and  the unknowns,
which are defined in the domain $\Omega\times (0,T_e)$, have the following form
$$
 b(t,x)  = \hat b(t,r)\frac{x}{r},\  \gamma(t,x) = \hat\gamma(t,r),\ S_0(x) = \hat S_0(r)
$$
and
$$
 u(t,x) = \hat u(t,r)\frac{x}{r},\ S(t,x) = \hat S(t,r),
$$
respectively, where $T_e$ is a positive constant which denotes the life-span of weak solutions,
$r=|x|$,  $ \Omega=\{x\in\cR^3\mid a< r <d\, \}$ for
 two positive constant $a,\ d$ satisfying $a<d$, and $\hat b,\ \hat\gamma,\ \hat S_0$ are given functions
 and $\hat u,\ \hat S$  are scalar functions to be determined, which
 depend only on $t,\ r$.    We write
$$
x=(x_i), \ u=(u_i),\ T=(T_{ij}),\ D=(D_{kl}^{ij}).
$$
Here and hereafter, $i,j,k,l=1,2,3$, and we assume that $D$ satisfies the properties of symmetry:
\begin{eqnarray}
 D_{kl}^{ij} = D^{kl}_{ij} = D_{lk}^{ij} = D_{kl}^{ji}.
 \label{AssumptionD1}
\end{eqnarray}
Moreover  we assume that $D$ satisfies
\begin{eqnarray}
 D_{ij}^{kl} &=& 0, \ {\rm if}\ k \not= j;
 \label{AssumptionD2}\\
 D_{ij}^{jl} &=& 0, \ {\rm if}\ i \not= l,\ {\rm for\ any\ fixed}\ j.\ {\rm Assume\ that}\ D_{ij}^{jl} \
 {\rm is\ independent\ of\ } j;
 \label{AssumptionD3}\\
 C_{il} &:=& D_{ij}^{jl}. \ {\rm Assume\ that}\ C_{ll} \
 {\rm is\ independent\ of\ } l\ {\rm and\ is\ equal\ to}\ \mu.\quad
 \label{Cil}\\
 E_{kl} &:=& \sum_{i,j=1}^3 D_{ij}^{kl}\bar\varepsilon_{ij}  =  0, \ {\rm if}\ k \not= l,\ {\rm and } \\
 & & E_{kk}\ {\rm
 is\ independent\ of \ } k\ {\rm and\ is\ equal\ to}\ \lambda.\ \ \
 \label{AssumptionEpsilon1}
\end{eqnarray}

\bigskip
Under these assumptions  equations \eq{eq1} -- \eq{eq2}  are reduced to
\begin{eqnarray}
 \frac{\partial^2}{\partial r^2}\hat u   + \frac2r\frac{\partial}{\partial r}\hat u
  - \frac2{r^2}\hat u   &=&  {\cal G},
 \label{sFinal1}\\
 \frac{\partial}{\partial t}\hat S
  + \left(- c\, \nu    \frac{\partial^2}{\partial r^2} \hat  S
  + {\cal F} \right)   |\frac{\partial}{\partial r}\hat S| &=& 0.
 \label{sFinal2}
\end{eqnarray}
Here ${\cal F}, {\cal G}$  are nonlinear functions defined  by
\begin{eqnarray}
 {\cal G} &=& {  {\cal G}(\frac{\partial}{\partial r}\hat S , \hat b) = \frac\lambda\mu \frac{\partial}{\partial r}\hat S + \frac{\hat b}{\mu} ,} \label{calG}\\
 {\cal F}_1 &=& {\cal F}_1  \Big(\hat u ,\frac{\partial}{\partial r}\hat u, \hat S,  \Big)\non\\[0.2cm]
 &=&   c \left({  -}\lambda \Big(\frac{\partial}{\partial r}\hat u + \frac{2}{r }\hat u \Big)
 + D\bar\varepsilon \cdot \bar\varepsilon \hat S + \hat\psi^\prime(\hat S)\right),\\
 {\cal F}  &=& {\cal F} \Big(\hat u ,\frac{\partial}{\partial r}\hat u, \hat S,  \frac{\partial}{\partial r}\hat S\Big)
 = {\cal F}_1 - \frac{2c\, \nu}{ r }\frac{\partial}{\partial r}\hat S.
 \label{sFinal2a}
\end{eqnarray}
The boundary and initial conditions become
\begin{eqnarray}
 \hat u(t,r)= \hat \gamma(t,r), & & \hat S(t,r)=0,\ (t,r) \in [0,T_e]\times \partial\Omega ,
 \label{sFinal3a}  \\
 && \hat S(0,r) = \hat S_0(r),\ r\in \Omega,
 \label{sFinal3}
\end{eqnarray}
where $ \hat \gamma(t,r)$ is defined by $\gamma(t,r) =  \hat \gamma(t,r)\frac{x}r$.

\medskip
\noindent{\bf Remark 1.} {\it One can easily find an example which meets the above assumptions:
 The media is isotropic and homogenous. These  assumptions  still lead to an elliptic-parabolic coupled system, thus
 the reduced system possesses the main difficulties in the proof
 of existence of weak solutions as in \cite{Alber06}.

}

\medskip
To what follows, except Section~2 in which we reduce the problem to one dimensional form, we
shall change  the independent variable $r$ to $x$, and drop the hat $\hat\ $ for
all quantities (except $\hat\psi$), namely,
$\hat u\to u,\ \hat b\to b$, etc. Denote $f_x = \frac{\partial}{\partial x} f$, $f_{xx} = \frac{\partial^2}{\partial x^2} f$, etc.
 The domain $\Omega$ is reduced to an interval: $\Omega=(a,d)$ is a
bounded open interval with constants $a<d$. We write $Q_{T_e}:=(0,T_e) \times
\Omega$, where $T_e$ is a positive constant, and define
$$
 (v,\varphi)_Z = \int_Z v(y)\varphi(y)\, dy\, ,
$$
for $Z = \Omega$ or $Z = Q_{T_e}$. If $v$ is a function defined on
$Q_{T_e}$ we denote the mapping $x \rightarrow v(t,x)$ by $v(t)$. If
no confusion is possible we sometimes drop the argument $t$ and
write $v = v(t)$. Since  equation \eq{sFinal1} is linear, the
inhomogeneous Dirichlet boundary condition for $\hat u$ can be reduced in the
standard way to the homogeneous condition. For simplicity we thus assume that
$$
\hat \gamma=0 .
$$

Then  with these simplifications, equations (\ref{sFinal1}) -- (\ref{sFinal2})   can be written in the form
\begin{eqnarray}
 u_{xx}  + \frac2x u_{x}
  - \frac2{x^2} u   &=& {\cal G},
 \label{Final1}\\
 \frac{\partial}{\partial t} S  + \left( {\cal F }
  - c \, \nu  S_{xx}   \right) | S_{x}| &=& 0.
 \label{Final2}
\end{eqnarray}
The boundary and initial conditions turn out to be
\begin{eqnarray}
  u(t,x)= 0,\ \ \   S(t,x) &=& 0,\ (t,x)\in [0,T_e]\times \partial\Omega ,
 \label{Final3a}  \\
  S(0,x) &=& S_0(x),\ x\in \Omega.
 \label{Final3}
\end{eqnarray}

\medskip
To define weak solutions of this initial-boundary value problem we
note that because of $\frac12(|y|y)^\prime=|y|$ equation
(\ref{Final2}) is equivalent to
\begin{eqnarray}
 \frac{\partial}{\partial t} S  - \frac{c\, \nu}{2} \left( S_{x}| S_{x}| \right)_{x} + {\cal F }  | S_{x}| = 0.
 \label{Final2Div}
\end{eqnarray}
%
\noindent{\bf Definition 1.1.} {\it Let $b\in L^\infty (0,T_e,L^2(\Omega))$,
$S_0\in L^\infty(\Omega)$. A function $(u, S)$ with
\begin{eqnarray}
 && u\in L^\infty(0,T_e;W^{1,\infty}_0(\Omega)) ,
 \label{property01} \\
 && S\in L^\infty(Q_{T_e})\cap L^\infty(0,T_e, H^1_0(\Omega)),
 \label{property02}
\end{eqnarray}
 is a weak solution to the problem (\ref{Final1}) -- \eq{Final3},
 if the equation (\ref{Final1}) is satisfied weakly and if for all $\varphi\in
C^\infty_0((-\infty,T_e)\times \Omega)$
\begin{equation}
 (S,\varphi_t)_{Q_{T_e}}-\frac{c\, \nu}2 ( |S_x|S_x,\,\varphi_x )_{Q_{T_e}} -
  ( {\cal F}   |S_x|,\, \varphi )_{Q_{T_e}} + (S_0,\varphi(0))_\Omega = 0.
 \label{definition}
\end{equation}
}
\\[1ex]
The main result of this article is
\begin{Theorem}\label{T1.1} To all $S_0\in H^1_0(\Omega)$ and $b \in
C(\overline{Q}_{T_e})$ with $b_t \in  C(\overline{Q}_{T_e})$ there
exists a weak solution $(u, S)$ of problem (\ref{Final1}) --
(\ref{Final3}), which in addition to (\ref{property01}) --
(\ref{definition}) satisfies
\begin{equation}
 S_t\in L^\frac43(Q_{T_e}),\quad  S_x\in  L^\frac83(0,T_e; L^\infty(\Omega)),
 \label{proper1}
\end{equation}
and
\begin{equation}
 (|S_x|S_x)_x\in L^\frac43(Q_{T_e}),\quad
 S_{xt}\in L^\frac43(0,T_e;W^{-1,\frac43}(\Omega)).
 \label{proper2}
\end{equation}

Consequently we find spherically symmetric solution
$(\hat u(t,r)\,\frac{x}{r},\hat S(t,r))$ to the original problem \eq{eq1} -- \eq{eq3}.

\end{Theorem}

\bigskip
The remaining sections are devoted to the proof of Theorem~\ref{T1.1}.
 The difficulties in the proof of existence of weak solutions to the one dimensional
 problem are due to the following features: The system is of
elliptic-parabolic type,  it consists of a linear second order elliptic equation
coupled with a nonlinear equation for the order parameter equation.  The nonlinearity of this nonlinear equation is
stronger than the one in \cite{Alber06}, where a one-dimensional initial boundary value problem
is investigated. Moreover, this equation is degenerate and the nonlinearity depends non-smoothly
on the gradient of unknown $S$. This can be judged easily from the fact that the coefficient
$\nu|S_x|$ of the highest order derivative $S_{xx}$ in the order parameter equation is not bounded
away from zero and that it is not differentiable with respect to $S_x$.

The rest of this article is organized as follows:  In Section~2, assuming that the domain $\Omega$, the
  elasticity tensor $D$ and the misfit stain tensor satisfy suitable conditions, and  that
  the solutions $(u,S)$ to problem \eq{eq1} -- \eq{eq3} and the initial and boundary data are spherically
  symmetric, we reduce the original problem to the one dimensional form.

Then to prove Theorem~\ref{T1.1} we first consider in Section~3 a modified
initial-boundary value problem which consists of (\ref{Final1}) and the equation
\begin{equation}
 S_t - c\, \nu |S_x|_\kappa S_{xx} +  {\cal F }\cdot   (
 |S_x|_\kappa - \kappa) =0 ,  \  x\in \Omega,\  t>0
 \label{eq2a}
\end{equation}
with a constant $\kappa>0$. Here we use the notation
\begin{equation}
 |p|_\kappa:= \sqrt{\kappa^2 + p^2} .
 \label{pfunction}
\end{equation}
Since (\ref{eq2a}) is a uniformly parabolic equation we can use a standard
theorem to conclude that the modified initial-boundary value problem has a
sufficiently smooth solution $(u^\kappa, S^\kappa)$. For this solution
we derive in Section~4 a-priori estimates that are uniform in  $\kappa$ for $\kappa\in (0,1]$. The
assumption $\kappa\in (0,1]$ is reasonable since we consider limits of approximate solutions for
 $\kappa\to 0$. We shall see that the selection of a function in the form \eq{pfunction} results in a simpler
 proof of the existence of weak solutions than that in \cite{Alber06}.

To select a subsequence converging to a solution for $\kappa\to 0$ we need a
compactness result. However, our a-priori estimates are not strong enough to
show that the sequence $S^\kappa_x$ is compact; instead, we can only show that
the sequence $\int_0^{S^\kappa_x}|y| dy = \frac12 S^\kappa_x |S^\kappa_x| $
has bounded derivatives, with respect to both $x$ and $t$, in some suitable spaces, and thus can be proved to be compact.
It turns out that this is enough to prove existence of a solution. For the compactness proof in
Section~5 we use the Aubin-Lions Lemma; since one of our a-priori estimates
for derivatives of the approximate solutions is only valid in
$L^1(0,T_e; H^{-2}(\Omega))$, we must use the generalized form of this lemma
given by Roub\'ic\v{e}k~\cite{Roubicek}, which is valid in $L^1$.

Despite we prove the existence of spherically symmetric solutions, the
 existence of weak solutions to the original problem \eq{eq1} -- \eq{eq3}
is still open. The method of the proof in this article and \cite{Alber06} (in which
 problem \eq{eq1} -- \eq{eq3} in one dimensional case is studied) is limited to one space dimension,
since for the a-priori estimates it is crucial that the term $|S_x|S_{xx}$ in (\ref{Final2}) can be
 written in the form $\frac12(|S_x|S_{x})_x$. In the higher dimensional case
the corresponding term $|\nabla_xS|\Delta_xS$ cannot be rewritten in this
way. Yet, we believe that these essentially one-dimensional existence results can also
be helpful in an existence proof for higher space dimensions.
%

\section{Reduction to one dimensional problem}

We shall prove in this section that under suitable assumptions, the
original  problem can be reduced to a one dimensional problem. We
now assume that the body force and  the unknowns, which are defined
in the domain $\Omega\times (0,T_e)$, have the following form
\begin{eqnarray}
 b(t,x)  = \hat b(t,r)\frac{x}{r},
  \label{form1}
\end{eqnarray}
and
\begin{eqnarray}
 u(t,x) = \hat u(t,r)\frac{x}{r},\ S(t,x) = \hat S(t,r),
  \label{form2}
\end{eqnarray}
respectively, where $T_e$ is a positive constant,  $r=|x|$,
$\Omega=\{x\in\cR^3\mid a< r <d\}$ for
 two positive constant $a,\ d$ satisfying $a<d$, and $\hat u,\ \hat S$
 are scalar functions to be determined, which
 depend only on $t,\ r$, and $  \hat b$ is a given function in $t,\ r$.

\begin{Theorem}\label{T0} Suppose that the tensors $D$ and $\bar\varepsilon$ satisfy
\eq{AssumptionD1} -- \eq{AssumptionEpsilon1}.
Then the following two statements are equivalent:

1. $(u,S)(t,x)$ of the form \eq{form2} is a classical solution to the problem \eq{eq1}
-- \eq{eq3} with $b$ chosen in \eq{form1},

2. $(\hat u,\hat S)(t,r)$ solves classically the problem \eq{sFinal1} -- \eq{sFinal3}.

\end{Theorem}

\noindent{\it Proof.}  To simplify notations, the Einstein summation
convention applies to the rest of this section: When an index
variable (e.g. $i,j,k,l$,
 but with an exception $r$ in this article, for instance, $\hat S_{,r} \frac{x_i}r$ in \eq{Reduction4}  does not
 mean that we take the sum for the index $r$) appears twice in a single term that is a product of two or more numbers,
 it implies that we are summing over all of its possible values. However we shall still
 use the symbol $\Sigma$ to avoid some possible confusion when in a single term, an index appears more than two times.
 Recalling
$$
 x=(x_i), \ u=(u_i),\ T=(T_{ij}),\ D=(D_{kl}^{ij}).
$$
where  $i,j,k,l=1,2,3$.

For partial derivatives, we denote for a function $\hat f=\hat f(t,r)$
$$
 \hat f_{,r} = \frac{\partial \hat f}{\partial r},\ \hat  f_{,rr} = \frac{\partial^2\hat  f}{\partial r^2}.
$$
An index $j$ (or $i, k,l,$) after a comma in subscript of a quantity (for example,
a function $f=f(t,x)$, a vector $u=u(t,x)$ and tensor $T=T(t,x)$, etc.)
 indicates the partial derivative with respect to $x_j$,  namely
$$
  f_{,j} = \frac{\partial f}{\partial x_j},\  u_{i,j} = \frac{\partial u_i}{\partial x_j},\  T_{ik,j} = \frac{\partial T_{ik}}{\partial x_j}, \cdots.
$$
Similar convention applies to multiple indices after a comma in subscript of a quantity.
 We can thus write
\begin{eqnarray}
 r_{,i} &=& \frac{x_i}{r},\quad \left(\frac{x_i}r\right)_{,j}= \frac{\delta_{ij}r^2 - x_ix_j}{r^3},
 \label{Reduction1}\\
 u_{i,j} &=& \left(\hat u  \frac{x_i}{r}\right)_{,j} = \hat u_{,r}\frac{x_ix_j}{r^2} + \hat u \frac{\delta_{ij}r^2 - x_ix_j}{r^3},
 \label{Reduction2}\\
 u_{i,jk} &=&  \hat u_{,rr}\frac{x_ix_jx_k}{r^3} + \hat u_{,r} \frac{r^2\left((x_ix_j)_{,k} + \delta_{ij}x_k \right) - 3x_ix_jx_k}{r^4} \non\\
  & & + \hat u\,  \frac{r^2\left( -\delta_{ij} x_k - (x_ix_j)_{,k} \right) + 3   x_ix_j  x_{k} }{r^5},
 \label{Reduction3} \\
 S_{,i} &=& \hat S_{,r} \frac{x_i}{r}.
  \label{Reduction4}
\end{eqnarray}
Here $\delta_{ij}$ is the Kronecker delta.

Hence, the first two equations \eq{eq1} -- \eq{eq1a} can be rewritten as
\begin{eqnarray}
\hat b\, \frac{x_l}{r} &=&  \frac12 D_{ij}^{kl} u_{j,ik}   + \frac12 D_{ij}^{kl}  u_{i,jk} - D_{ij}^{kl}\bar\varepsilon_{ij} S_{,k}\non\\
 &=&  \hat u_{,rr} D_{ij}^{kl}  \frac{x_ix_jx_k}{r^3} + \hat u_{,r} \left( D_{ij}^{kl}  \frac{ (x_ix_j)_{,k}}{r^2} + D_{ij}^{kl}  \frac{ \delta_{ij}x_k   }{r^2} - D_{ij}^{kl}  \frac{ 3x_ix_jx_k}{r^4} \right) \non\\
  & & +   \frac{\hat u}{r^5} \left(-  r^2 \left( D_{ij}^{kl} \delta_{ij} x_k + D_{ij}^{kl} (x_ix_j)_{,k}
 \right)  +  3 D_{ij}^{kl}   x_ix_jx_{k}  \right) - D_{ij}^{kl}\bar\varepsilon_{ij}\hat S_{,r}\frac{ x_k}{r }.
 \label{Reduction5}
\end{eqnarray}

\medskip
From now on we take one $l$ from $\{1,2,3\}$.
 Invoking assumptions \eq{AssumptionD2} -- \eq{AssumptionEpsilon1} we obtain
\begin{eqnarray}
 D_{ij}^{kl}  \frac{x_ix_jx_k}{r^3} = \sum_{j=1}^3 D_{ij}^{jl}  \frac{x_ix_jx_j}{r^3}
 = C_{il}  \frac{x_i}{r } = C_{ll}  \frac{x_l}{r } = \mu \frac{x_l}{r }.
 \label{Reduction5a}
\end{eqnarray}
Since $(x_ix_j)_{,k} =  \delta_{ik} x_j +  x_i \delta_{jk}$,  one has
\begin{eqnarray}
 && D_{ij}^{kl}  \frac{ (x_ix_j)_{,k}}{r^2} + D_{ij}^{kl}  \frac{ \delta_{ij}x_k   }{r^2} - D_{ij}^{kl}  \frac{ 3x_ix_jx_k}{r^4}
   \non\\
 &=& \sum_{j=1}^3\left( D_{ij}^{jl}  \frac{ \delta_{ij} x_j +  x_i \delta_{jj}}{r^2} + D_{ij}^{jl}  \frac{ \delta_{ij}x_j   }{r^2} - D_{ij}^{jl}  \frac{ 3x_ix_jx_j}{r^4}
  \right) \non\\
 &=&\sum_{i=1}^3 C_{il} \frac{r^2\left( 3 x_i  + 2 x_i \right) - 3x_i r^2}{r^4} = C_{il}\frac{2r^2x_i }{r^4}\non\\
 &=&   2 \mu\frac{ x_l }{r^2},
 \label{Reduction5b}
\end{eqnarray}
and
\begin{eqnarray}
 && -\frac{  1}{r^3}  \left( D_{ij}^{kl} \delta_{ij} x_k + D_{ij}^{kl} (x_ix_j)_{,k}   \right)
  +   D_{ij}^{kl} \frac{ 3 x_ix_jx_{k} }{r^5}\non\\
 &=& \sum_{j=1}^3 D_{ij}^{jl}   \frac{ -r^2\left( \delta_{ij} x_j + (x_ix_j)_{,j} \right) + 3 x_ix_j  x_{j} }{r^5}\non\\
 &=& \sum_{j=1}^3 D_{ij}^{jl}   \frac{ -r^2\left( 2\delta_{ij} x_j + x_i\delta_{jj} \right) + 3 x_ix_j  x_{j} }{r^5}\non\\
 &=& \sum_{i=1}^3 C_{il}   \frac{-r^2(  2 x_i + 3x_i ) {  + } 3  r^2 x_{i} }{r^5}
 = C_{il}\frac{-2 x_i  }{r^3}\non\\
 &=& -2\mu \frac{ x_l  }{r^3}.
 \label{Reduction5c}
\end{eqnarray}
Using \eq{Reduction5a} -- \eq{Reduction5c}, we are in a position to rewrite equation \eq{Reduction5} as
\begin{eqnarray}
 \hat b \frac{x_l}{r} &=& \frac{x_l}{r}\left( \mu\Big( \hat u_{,rr} +
 \frac2r\hat u_{,r} - \frac2{r^2}\hat u \Big) - \lambda
 \hat S_{,r}\right).
 \label{Reduction6}
\end{eqnarray}
This holds, for $r\ge a>0$, if and only if the following equation is satisfied
\begin{eqnarray}
 \mu \Big( \hat u_{,rr}   + \frac2r\hat u_{,r}
  - \frac2{r^2}\hat u \Big) - \lambda \hat S_{,r}  &=&  \hat b,
 \label{ReductionFinal1}
\end{eqnarray}
which is just \eq{sFinal1}.

\medskip
Next, to deal with the order parameter equation we make use of the following formula
$$
 \psi_S(\varepsilon, S)=-T\cdot \bar\varepsilon + \hat\psi^\prime(S)
 = -D\varepsilon(\nabla u)\cdot \bar\varepsilon + D\bar\varepsilon \cdot \bar\varepsilon S + \hat\psi^\prime(S) .
$$
We evaluate $D\varepsilon(\nabla u)\cdot \bar\varepsilon $. Invoking \eq{AssumptionEpsilon1},
\begin{eqnarray}
 D\varepsilon(\nabla u)\cdot \bar\varepsilon &=& \frac12 D_{ij}^{kl} u_{i,j}  \bar\varepsilon_{kl}
  + \frac12 D_{ij}^{kl}  u_{j,i} \bar\varepsilon_{kl}\non\\
 &=& \hat u_{,r} D_{ij}^{kl}  \frac{x_ix_j}{r^2}\bar\varepsilon_{kl} +
 \hat u \left( D_{ij}^{kl}\frac{\delta_{ij} }{r }\bar\varepsilon_{kl}
 - D_{ij}^{kl}\frac{  x_ix_j}{r^3}  \bar\varepsilon_{kl} \right)\non\\
 &=& \hat u_{,r} E_{ij} \frac{x_ix_j}{r^2} +
 \hat u \left( E_{ij} \frac{\delta_{ij} }{r }
 - E_{ij} \frac{  x_ix_j}{r^3} \right)\non\\
 &=&  \sum_{i=1}^3 E_{ii} \left(\hat u_{,r}\frac{x_ix_i}{r^2} + \hat u \frac{\delta_{ii}r^2 - x_ix_i}{r^3} \right)\non\\
 &=&  \lambda  \left(\hat u_{,r} + \frac{2}{r }\hat u \right),
 \label{Reduction6a}
\end{eqnarray}
where $\lambda$ is a constant as in \eq{AssumptionEpsilon1}. Thus \eq{eq2} turns out to be
\begin{eqnarray}
 \hat S_t  +  (- c\, \nu \hat S_{,rr} +   {\cal F }) |\hat S_{,r}|= 0,
 \label{ReductionFinal2}
\end{eqnarray}
where ${\cal F } $ is the same function as  in \eq{sFinal2a}. Thus we obtain \eq{sFinal2}.

To finish the reduction of the problem, we write the initial boundary
conditions in the following form: $\gamma(t,x) = \hat
\gamma(t,r)\frac{x}{r}$, $S_0(x)=\hat S_0(r)$. Thus we obtain the
one dimensional problem and the proof of Theorem~\ref{T0} is complete.

\section{Existence of solutions to the modified problem}
In this section, we study the modified initial-boundary value problem and
show that it has a H\"older continuous classical solution, consequently we
construct approximate solutions whose limit is a solution to the original
problem \eq{Final1} -- \eq{Final3}. To formulate this problem, let
$
\chi\in C_0^\infty(\cR,[0,\infty))
$
satisfy
$
\int_{-\infty}^{\infty} \chi(t)dt=1.
$
For $\kappa>0$, we set
$$
 \chi_\kappa(t):=\frac1\kappa\chi\left(\frac{t}{\kappa} \right),
$$
and for $S\in L^\infty(Q_{T_e},\cR)$ we define
\begin{eqnarray}
 (\chi_\kappa*S)(t,x)=\int_{0}^{T_e} \chi_\kappa(t-s)S(s,x) ds.
 \label{convolution}
\end{eqnarray}

The modified initial-boundary value problem consists of the equations
\begin{eqnarray}
 u_{xx}   + \frac2x u_{x} - \frac2{x^2} u  &=&  {\cal G}( (\chi_\kappa*S) _{ x}, b) ,
 \label{m2.1}\\
 S_t -  c\, \nu|S_x|_\kappa S_{xx} &=& -  {\cal F}\cdot (|S_x|_\kappa- \kappa),
 \label{m2.2}
\end{eqnarray}
which must hold in $Q_{T_e}$, and of the boundary and initial conditions
\begin{eqnarray}
 u(t,x) &=& 0,\ S(t,x) = 0,\  (t,x)\in (0,T_e)\times \partial\Omega,
 \label{m2.4}\\
 S(0,x) &=& S_0(x), \quad x\in \bar\Omega.
 \label{m2.6}
\end{eqnarray}

\medskip
Now we want to rewrite the system as an equation with a nonlocal term.
 Applying the Sturm-Liouville theory for ordinary differential equations of the form
$$
\frac{d}{dx}\left(p(x) y_{x}(x)\right) +  q(x)y(x) = 0,
$$
with suitable boundary conditions at $x=a, d$, we first solve $u$ in terms of $S_x$ and $b$. To this end,
we rewrite, by multiplying it by $x^2$, \eq{ReductionFinal1} as
\begin{eqnarray}
 L[ u] &:=& \frac{d}{dx}\left(  p(x)  u_{x}\right) + q(x) u ,
  \label{m2.1a} \\
  L[ u](x)  &=&  {x^2}{\cal G}( (\chi_\kappa*S) _{ x}, b) ,
 \label{m2.1b}
 \end{eqnarray}
and the boundary conditions are chosen as $  u (t,x) =0$ at $x=a,d$. Here
$$
 p(x)= x^2, \ q(x)= -  2 .
$$

Consider the eigen-problem $L[ u]=\sigma  u$ with $\sigma=0$ and with $ u (t,x) =0$ at $x=a,d$.
It is easy to show from \eq{m2.1b} that
$$
0 = \int_a^d L[  u] \cdot u dx = \int_a^d \left( \frac{d}{dx}\left( p(x) u_{x}\right) u
  + q(x)  u^2 \right)dx   = - \int_a^d \left( x^2  u_{x}^2
  +  2 u^2 \right)dx ,
$$
 whence $  u\equiv 0$, and $ 0$ is not an eigenvalue of this operator.
 One asserts that for any fixed $t\in [0,T_e]$, there
exists a unique solution $ u$ to \eq{m2.1}, which can be represented by
\begin{eqnarray}
 u (t,x) =   \int_a^d G(x,y)\left(\frac{y^2}\mu b(t,y) + \frac\lambda\mu  y^2 (\chi_\kappa* S(t,y) )_{y}\right) dy.
 \label{m2.1c}
\end{eqnarray}
Here $G(x,y)$ is the Green function, associated with the operator $L$, such that \\
1. $G(x,y)$ is continuous in $x$ and $y$;\\
2. For $x\not=y$, $L[G(x,y)]=0$;\\
3.   $ G(x,\cdot) =0$ at $x=a,d$;\\
4. Derivative jump: $G'(y_{+0},y) - G'(y_{-0},y) = \frac1{p(y)}$;\\
5. Symmetry: $G(x,y) = G(y,x)$.

Recalling the boundary condition \eq{m2.4} and integrating by parts we infer from \eq{m2.1c} that
\begin{eqnarray}
  u (t,{  x}) & = & \frac1\mu \int_a^d G(x,y) y^2  b(t,y) dy - \frac\lambda\mu  \int_a^d \left(G(x,y) y^2\right)_{y}  \chi_\kappa* S (t,y)  dy\non\\
  & = & \frac1\mu \int_a^d G(x,y)y^2 b(t,y) dy - \frac\lambda\mu   \int_{ \{ x\not= y\} } 2G(x,y) y \chi_\kappa* S(t,y) dy\non\\
  & &   -\frac\lambda\mu   \int_{ \{x\not= y\} }   G(x,y)_{y} y^2 \chi_\kappa* S(t,y) dy  .
  \label{m2.1d}
\end{eqnarray}
Thus $u (t,x)$ depends linearly on $S$ and a nonlocal term of $S$.

\medskip
To formulate an existence theorem for this problem we need some
function spaces: For nonnegative integers $m,n$ and a real number
$\alpha\in (0,1)$ we denote by $C^{m+\alpha}(\overline{\Omega})$  the space of
$m-$times differentiable functions on $\overline{\Omega}$, whose $m-$th
derivative is H\"older continuous with exponent $\alpha$.  The space
$C^{\alpha,\alpha/2}(\overline{Q}_{T_e})$ consists of all functions on
$\overline{Q}_{T_e}$, which are H\"older continuous in the parabolic distance
$$
 {\rm d}((t,x),(s,y)):=\sqrt{|t-s|+|x-y|^2}.
$$
$C^{m,n}(\overline{Q}_{T_e})$ and
$C^{m+\alpha,n+\alpha/2}(\overline{Q}_{T_e})$, respectively, are the spaces of
functions, whose $x$--derivatives up to order $m$ and $t$--derivatives up to
order $n$ belong to $C(\overline{Q}_{T_e})$ or to
$C^{\alpha,\alpha/2}(\overline{Q}_{T_e})$, respectively.
%
\begin{Theorem}\label{T2.1} Let $\nu,\kappa>0$, $T_e > 0$, suppose
that the function $b \in C(\overline{Q}_{T_e})$ has the derivative
$b_t \in C(\overline{Q}_{T_e})$ and that the initial data $S_0\in
C^{2+\alpha}(\overline{\Omega})$ satisfy $S_0|_{\partial\Omega} =
S_{0,x}|_{\partial\Omega} = S_{0,xx}|_{\partial\Omega} = 0$. Then
there is a  solution
$$
 (u, S) \in C^{2,1}(\overline{Q}_{T_e})  \times
 C^{2+\alpha,1+\alpha/2}(\overline{Q}_{T_e})
$$
of the modified initial-boundary value problem (\ref{m2.1}) --
(\ref{m2.6}). This solution satisfies $S_{tx}\in L^2(Q_{T_e})$ and
\begin{eqnarray}
 \max_{\overline{Q}_{T_e}}|S|\le \max_{\overline{\Omega}}|S_0|.
 \label{m2.7}
\end{eqnarray}

\end{Theorem}
\vskip0.2cm
{\it Proof.} Making use of \eq{m2.1d}, we  rewrite the system \eq{m2.1} -- \eq{m2.2} as a single equation
\begin{eqnarray}
 S_t &=& a_1(S_x)S_{xx}\non\\
 && + a_2\left(t,x,S,S_x,\tilde S,
  \int_{ \{ x\not= y\} }  G(x,y) y \tilde S (t,y) dy,  \int_{ \{x\not= y\} }   G(x,y)_{y} y^2 \tilde S(t,y) dy \right)\qquad
\label{m2.9}
\end{eqnarray}
in $Q_{T_e}$\,, where $\tilde S = \chi_\kappa*S$,
\begin{eqnarray*}
a_1(p)=c\, \nu|p|_\kappa
\end{eqnarray*}
and
\begin{eqnarray*}
a_2(t,x,S,p,r,s) = c\, {  \overline{\cal F} } (t,x,S,p,r,s_1,s_2)
(|p|_\kappa - \kappa).
\end{eqnarray*}
Here ${\overline{\cal F} }$ is obtained by using formula \eq{m2.1d}
and inserting $u,u_x$ into the formula of ${\cal F}$.

Equation (\ref{m2.9}) is quasilinear, uniformly   parabolic equation,
 which  contains nonlocal terms.
Then  we can apply, with a little modification,
\cite[Theorem~2.9, p.23]{Ladyzenskaya} to (\ref{m2.9}) to prove the existence of classical solution $S^\kappa$, and conclude
that the estimate (\ref{m2.7}) holds by applying the maximum principle to
\eq{m2.9}. We refer the reader to the paper \cite{Alber06} for the details. Thus we complete the
 proof of Theorem~\ref{T2.1}.
%
%
\section{A priori estimates}
This section is devoted to the derivation of a-priori estimates for
solutions of the modified problem, which are uniform with respect to
$\kappa\in (0,1]$. We remark that the estimates in Lemma~\ref{L3.1}
and Corollary~3.1, though stated in the one-dimensional case, can be
generalized to higher space dimensions.

In what follows we assume that
\begin{equation}
 0<\kappa\le 1,
 \label{assumkappa}
\end{equation}
since we consider the limit $\kappa\to 0$. The $L^2(\Omega)$-norm is denoted
by $\|\cdot\|$, and the letter $C$ stands for varies positive constants
independent of $\kappa$.
We start by constructing a family of approximate solutions to the
modified problem. To this end let $T_e$ be a fixed positive
number and choose for every $\kappa$ a function $S_0^\kappa\in
C_0^\infty(\Omega)$ such that
\begin{equation}
 \|S_0^\kappa-S_0\|_{H^1_0(\Omega)}\to 0,\quad \kappa\to 0,
 \label{approximate}
\end{equation}
where $S_0\in H^1_0(\Omega)$ are the initial data given in Theorem~1.1. We insert for $S_0$ in (\ref{m2.6})
 the function $S_0^\kappa$ and
choose for $b$ in (\ref{m2.1}) the function given in
Theorem~1.1. These functions satisfy the assumptions of Theorem~\ref{T2.1},
hence there exists a solution $(u^\kappa, S^\kappa)$ of the
modified problem (\ref{m2.1}) -- (\ref{m2.6}), which exists in
$Q_{T_e}$. The inequality (\ref{m2.7}) and Sobolev's embedding theorem
yield for this solution
\begin{equation}\label{approximateA}
 \sup_{0<\kappa\le 1}\|S^\kappa \|_{L^\infty(Q_{T_e})} \leq
 \sup_{0<\kappa\le 1} \|S_0^\kappa\|_{L^\infty(\Omega)} \le C.
\end{equation}

Remembering the formula (\ref{m2.1d}) and assumptions of $b$, we
show easily that $u^\kappa$ belongs to $C^{1,1}(\bar Q_{T_e})$, and
conclude from (\ref{m2.1}) that also
$\|u^\kappa_x\|_{L^\infty(Q_{T_e})}\le C$, and invoke the definition
of ${\cal F}_1$ to get
\begin{equation}\label{5max}
 \max_{\overline{Q}_{T_e}} |{\cal F}_1(u^\kappa,u^\kappa_x,S^\kappa)| \le C .
\end{equation}
With the help of this  estimate, we can evaluate derivatives of $S^\kappa$.
%
\begin{Lemma} \label{L3.1}
There holds for any $t\in [0,T_e]$
 \begin{equation}
 \|S^\kappa_x(t)\|^2 + c\,\nu\int_0^t\int_\Omega |S^\kappa_x|_\kappa
 |S^\kappa_{xx}|^2dxd\tau \le C.
 \label{5}
 \end{equation}

\end{Lemma}

\noindent{\it Proof.} From the assertion that $S^\kappa_{tx} \in L^2(Q_{T_e})$ in Theorem~\ref{T2.1}, it follows that
 there holds for almost all $t$
$$
 \frac12\frac{d}{dt}\|S^\kappa_{x}(t)\|^2 = \int_\Omega S^\kappa_x(t)
 S^\kappa_{xt}(t) dx.
$$
Making use of this relation and (\ref{5max}) we obtain by multiplication of (\ref{m2.2})
by $-S^\kappa_{xx}$ and integration by parts with respect to $x$,
where we take the boundary condition (\ref{m2.4}) into account, that
for almost all $t$
\begin{eqnarray}
 &&\frac12\frac{d}{dt}\|S^\kappa_{x}\|^2+
 \int_\Omega c\nu|S^\kappa_{x}|_\kappa |S^\kappa_{xx}|^2dx
 =\int_\Omega  {\cal F } ( |S_x|_\kappa - \kappa) S^\kappa_{xx} dx \non
 \\[1ex]
 &\le& C\int_\Omega (1+|S^\kappa_{x}|)(|S^\kappa_{x}|_\kappa+\kappa) |S^\kappa_{xx}| dx\non\\
 &=& C\left( \int_\Omega |S^\kappa_{x}|_\kappa |S^\kappa_{xx}| dx + \int_\Omega \kappa |S^\kappa_{xx}| dx + \int_\Omega |S^\kappa_{x}|\,|S^\kappa_{x}|_\kappa  |S^\kappa_{xx}| dx + \int_\Omega \kappa |S^\kappa_{x}|\, |S^\kappa_{xx}| dx \right)\non\\
 &=& C(I_1+I_2+ I_3+I_4)
 \label{4}
\end{eqnarray}
Now we estimate $I_i$ ($i=1,2,3,4$). For $I_1$, we have
\begin{eqnarray}
 I_1 &\le& C\int_\Omega |S^\kappa_{x}|_\kappa^\frac12
  (|S^\kappa_{x}|_\kappa^\frac12| S^\kappa_{xx}|) dx \non\\[1ex]
 & \le & \frac{c\nu }8\int_\Omega|S^\kappa_{x}|_\kappa|S^\kappa_{xx}|^2 dx
 +  C_{\nu} \int_\Omega (|S^\kappa_x|_\kappa)^2 dx {  + C},
 \label{I1}
\end{eqnarray}
where we denote $C_{\nu }$   a constant depending on $ \nu $.
By definition, there holds
\begin{eqnarray}
|S^\kappa_{x}|_\kappa\ge \kappa.
 \label{lowerBound}
\end{eqnarray}
Thus we can use the second term on the left hand side of \eq{4} to absorb $I_2$. By the Cauchy-Schwarz and Young inequalities,
\begin{eqnarray}
 I_2 &\le& C \int_\Omega \kappa^\frac12 (\kappa^\frac12 |S^\kappa_{xx}|) dx
 \le C\left(\int_\Omega \kappa  dx\right)^\frac12 \left( \int_\Omega \kappa |S^\kappa_{xx}|^2 dx\right)^\frac12\non\\
 & \le & \frac{c\nu \kappa }8 \int_\Omega  |S^\kappa_{xx}|^2 dx + {  C_{\nu}}.
 \label{I2}
\end{eqnarray}
Here we have used the fact that $0<a\le x\le d$ which implies the term $\frac{1}{x}$ contained in ${\cal F }$
is uniformly bounded from below and above, and $C_{\nu,\kappa}$ is a constant depending on $ {\nu,\kappa}$.
 Moreover, $I_3$ is evaluated by
\begin{eqnarray}
 I_4 &\le& C \int_\Omega \kappa |S^\kappa_{x}|\, |S^\kappa_{xx}| dx
 \le  C \left(\int_\Omega \kappa |S^\kappa_{x}|^2 dx\right)^\frac12\left( \int_\Omega \kappa  |S^\kappa_{xx}|^2 dx\right)^\frac12.
 \label{I4}
\end{eqnarray}
Using the Nirenberg inequality in the following form
\begin{eqnarray}
 \|f_x\|\le C \|f_{xx}\|^\frac13\|f\|_{L^\infty(\Omega)}^\frac23 + C' \|f\|_{L^\infty(\Omega)},
 \label{Nirenberg1}
\end{eqnarray}
one infers from \eq{I4} that
\begin{eqnarray}
 I_4 &\le&  C  \kappa^\frac12 \left(\|S^\kappa_{xx}\|^\frac13\|S^\kappa\|_{L^\infty(\Omega)}^\frac23 + C' \|S^\kappa\|_{L^\infty(\Omega)}\right) \left( \int_\Omega \kappa  |S^\kappa_{xx}|^2 dx\right)^\frac12\non\\
 &\le&  C  \kappa ^\frac12 \left(\|S^\kappa_{xx}\|^\frac13  + 1 \right) \left( \int_\Omega \kappa  |S^\kappa_{xx}|^2 dx\right)^\frac12\non\\
 &=&  C  \kappa \left(\|S^\kappa_{xx}\|^{\frac13  + 1} +   \|S^\kappa_{xx}\| \right) \non\\
  &\le& \frac{c\nu \kappa }8 \int_\Omega |S^\kappa_{xx}|^2 dx + {  C_{\nu}}\, ,
 \label{I4a}
\end{eqnarray}
where we used the Young inequality of the form $ab\le \delta a^\frac43 + C_\delta b^4 $ for non-negative real numbers  $a,b$. Here $C_\delta$ is a
constant depending on $ \delta$.

$I_3$ is the most difficult term to deal with. Again by the Cauchy-Schwarz inequality one gets
\begin{eqnarray}
 I_3 &=&  C \int_\Omega (|S^\kappa_{x}| |S^\kappa_{x}|_\kappa^\frac12)( |S^\kappa_{x}|_\kappa^\frac12 |S^\kappa_{xx}|) dx \non\\
 &\le&  C \left(\int_\Omega  |S^\kappa_{x}|^2 |S^\kappa_{x}|_\kappa   dx\right)^\frac12\left(  \int_\Omega |S^\kappa_{x}|_\kappa |S^\kappa_{xx}|^2 dx\right)^\frac12 \non\\
 &=&  C J^\frac12 \cdot \left(  \int_\Omega |S^\kappa_{x}|_\kappa |S^\kappa_{xx}|^2 dx\right)^\frac12
 \label{I3a}
\end{eqnarray}
To deal with  $J$, we recall the boundary conditions for $S^\kappa$ and rewrite by integration by parts
\begin{eqnarray}
 J & = & \int_\Omega  (S^\kappa_{x}) ^2 |S^\kappa_{x}|_\kappa   dx
 = \int_\Omega    S^\kappa_{x} (S^\kappa_{x} |S^\kappa_{x}|_\kappa )   dx \non\\
 & = & - \int_\Omega    S^\kappa  S^\kappa_{xx} \left(|S^\kappa_{x}|_\kappa + S^\kappa_{x} (|y|_\kappa)'|_{y=S^\kappa_{x}} \right)   dx\non\\
 & = & -\int_\Omega     S^\kappa S^\kappa_{xx}  \left( |S^\kappa_{x}|_\kappa + \frac{(S^\kappa_{x})^2}{ |S^\kappa_{x}|_\kappa }\right)  dx\non\\
 \label{I3b}
\end{eqnarray}
Applying estimate \eq{approximateA} and invoking the definition of
$|y|_\kappa$, from \eq{I3b} one obtains
\begin{eqnarray}
 J & \le &   C\int_\Omega    2 |S^\kappa_{x}|_\kappa| S^\kappa_{xx} | dx \non\\
 & \le & C\left(\int_\Omega |S^\kappa_{x}|_\kappa | S^\kappa_{xx} |^2 dx \right)^\frac12
   \left(\int_\Omega |S^\kappa_{x}|_\kappa    dx\right)^\frac12 \non\\
 & \le & C\left(\int_\Omega |S^\kappa_{x}|_\kappa | S^\kappa_{xx} |^2dx \right)^\frac12
   \left(\|S^\kappa_{x}\|+1 \right).
 \label{I3c}
\end{eqnarray}
Therefore, \eq{I3a} becomes
\begin{eqnarray}
 I_3 &\le& C  (\|S^\kappa_{x}\|+1)^\frac12 \left( \int_\Omega |S^\kappa_{x}|_\kappa |S^\kappa_{xx}|^2 dx\right)^{ \frac14 + \frac12 }\non\\
  &\le&   \frac{c\nu }8 \int_\Omega |S^\kappa_{x}|_\kappa |S^\kappa_{xx}|^2 dx  + C_\nu  (\|S^\kappa_{x}\|^2 + 1).
  \label{I3}
\end{eqnarray}
Here we have again used the Young inequality: $ab\le \delta a^\frac43 + C_\delta b^4 $.

Combining estimates \eq{I1} -- \eq{I3}, subtracting the terms
$\frac{c\nu}4\int_\Omega|S^\kappa_{x}|_\kappa|S^\kappa_{xx}|^2 dx$ and
$\frac{c\nu\kappa}4\int_\Omega |S^\kappa_{xx}|^2 dx$ on both sides of inequality \eq{4},
splitting the second term on the left hand side of \eq{4} into two equal terms, and  recalling the property \eq{lowerBound}, we derive
\begin{eqnarray}
 \frac12\frac{d}{dt}\|S^\kappa_{x}\|^2 +
 \frac{c\nu}4\int_\Omega |S^\kappa_{x}|_\kappa |S^\kappa_{xx}|^2dx + \frac{c\nu\kappa}4\int_\Omega |S^\kappa_{xx}|^2dx
  \le  C  \|S^\kappa_{x}\|^2 +  {  C_{\nu}}\,.
 \label{4a}
\end{eqnarray}
Then using the Gronwall inequality, one gets
(\ref{5}), noting also (\ref{approximate}). And the proof of this lemma is thus complete.

\medskip
Furthermore, we obtain
%
%
\begin{Corollary}\label{C3.2} There holds for any $t\in
[0,T_e]$
\begin{eqnarray}
 \int_0^t\int_\Omega\left(|S^\kappa_x|_\kappa|S^\kappa_{xx} |  \right)^\frac43dxd\tau &\le& C,
 \label{5d}\\
 \int_0^t\int_\Omega\left(|S^\kappa_x S^\kappa_{xx} |  \right)^\frac43dxd\tau &\le& C,
 \label{5w}\\
 \int_0^t\left\|\int_0^{S^\kappa_x}|y|_\kappa dy\right\|_{W^{1,\frac43}(\Omega)}^\frac43 d\tau &\le& C,
 \label{5g}\\
 \int_0^t \left\|\int_0^{S^\kappa_x}|y|_\kappa dy\right\|_{L^\infty(\Omega)}^\frac43 d\tau &\le& C,
 \label{5e}\\
 \|\, |S^\kappa_x| S^\kappa_x \|_{L^\frac43(0,T_e;L^\infty(\Omega))} &\le& C,
 \label{5h}\\
 \int_0^t \left\|S^\kappa_x\right\|_{L^\infty(\Omega)}^\frac83 d\tau &\le& C.
 \label{5f}
\end{eqnarray}
\end{Corollary}
{\it Proof.} For some $2 > p \ge 1$ we choose $q,\ q'$ such that
$$
 q=\frac2p, \quad \frac1q + \frac{1}{q^\prime}=1.
$$
By H\"older's inequality,  we have
\begin{eqnarray}
& & \int_0^t\int_\Omega\left(|S^\kappa_{x}|_\kappa|S^\kappa_{xx}|
  \right)^pdx d\tau \non\\
&= & \int_0^t\int_\Omega \left(|S^\kappa_{x}|_\kappa\right)^\frac{p}{2}
  \left( \left(|S^\kappa_{x}|_\kappa\right)^\frac{p}{2}
  |S^\kappa_{xx} |^p\right) dx d\tau \non \\
&\le&  \left(\int_0^t\int_\Omega
  \left(|S^\kappa_{x}|_\kappa \right)^\frac{pq^\prime}{2}  dxd\tau
  \right)^\frac{1}{q^\prime}\left(\int_0^t\int_\Omega
  \left(|S^\kappa_{x}|_\kappa \right)^\frac{pq}{2}
  |S^\kappa_{xx}|^{pq}  dxd\tau \right)^\frac{1}{q} \non\\
&\le &  \left(\int_0^t \int_\Omega
  \left(|S^\kappa_{x}|_\kappa\right)^\frac{p}{2-p}
  dx d\tau\right)^\frac{2-p}{2}
  \left(\int_0^t \int_\Omega |S^\kappa_{x}|_\kappa|S^\kappa_{xx} |^2
  dx d\tau \right)^\frac{p}{2}.
\label{6}
\end{eqnarray}
Estimate (\ref{5}) implies that if $p$ satisfies $\frac{p}{2-p}\le 2$ (i.e. $p\le
 \frac43$)  then the right hand side of (\ref{6}) is bounded. This
yields  estimate \eq{5d}.
 %
{Then \eq{5w} follows from \eq{5d} and the estimate \eq{5}, and the fact that $|S_x^\kappa|< |S_x^\kappa|_\kappa$.}

Next we are going to prove  \eq{5g}. Writing
\begin{eqnarray}
 |S^\kappa_x|_\kappa S^\kappa_{xx}= \left(\int_0^{S^\kappa_x}|y|_\kappa dy\right)_x,
 \label{sxxsx}
\end{eqnarray}
 and invoking that the primitive of $|y|_\kappa$ is equal to
$$
 \frac12\left(y\sqrt{y^2+ \kappa^2} + \kappa^2\log\Big(y + \sqrt{y^2+ \kappa^2} \Big)  \right),
$$
which,  thanks to $\log x\le  x -1 $ for all $x>0$, is bounded by $ C(y^2 + 1)$,  we then show easily that
$$
 \int_\Omega \int_0^{S^\kappa_x}|y|_\kappa dy dx\le C\int_\Omega (|S^\kappa_x|^2 +1)dx\le C.
$$
To apply the Poincar\'e  inequality of the form
$$
 \|f-\bar f\|_{L^p(\Omega)}\le C \|f_x\|_{L^p(\Omega)}
$$
where  $\bar f:=\frac1{|\Omega|}\int_\Omega f(x)dx$, we choose
$$
 p=\frac43,\quad f = \int_0^{S^\kappa_x}|y|_\kappa dy ,
$$
and obtain
\begin{eqnarray}
 &&\int_0^t \left\|\int_0^{S^\kappa_x}|y|_\kappa dy  \right\|^\frac43_{L^{\frac43} (\Omega) }d\tau  \non\\
 &\le& C\int_0^t \left\| \left(\int_0^{S^\kappa_x}|y|_\kappa dy \right)_x \right\|^\frac43_{L^{\frac43} (\Omega) }d\tau
 + C \int_0^t \left\|\, {\overline{ \int_0^{S^\kappa_x} |y|_\kappa   dy} }\, \right\|^\frac43_{L^{\frac43} (\Omega) } d\tau \non\\
 &\le& C\int_0^t \left\| \, |S^\kappa_x |_\kappa S^\kappa_{xx} \right\|^\frac43_{L^{\frac43} (\Omega) } d\tau  + C \int_0^t  1\, d\tau,
 \label{sxIntegral}
\end{eqnarray}
which implies, by \eq{5d}, that
\begin{eqnarray}
 \int_0^t \left\|\int_0^{S^\kappa_x}|y|_\kappa dy  \right\|^\frac43_{L^{\frac43} (\Omega) }d\tau
 &\le& C.
  \label{sxIntegral2}
\end{eqnarray}
Hence \eq{5g} follows, and we get $\int_0^{S^\kappa_x}|y|_\kappa dy\in L^\frac43(0,T_e;W^{1,\frac43}(\Omega))$.
Making use of the Sobolev embedding theorem, we get \eq{5e}.

It remains to prove estimate (\ref{5h}), since (\ref{5f})  is
equivalent to (\ref{5h}).
 We rewrite $\int_0^{S^\kappa_x}|y|_\kappa dy$ as
\begin{eqnarray}
\int_0^{S^\kappa_x}|y|_\kappa dy
 &= &\int_0^{S^\kappa_x}|y| dy+ \int_0^{S^\kappa_x}(|y|_\kappa-|y|) dy\non\\
 &= &\left.\frac12|y| y\right|_0^{S^\kappa_x}
+ \int_0^{S^\kappa_x}\frac{\kappa^2}{|y|_\kappa+|y|} dy\non\\
 &= &\frac12 |S^\kappa_x| S^\kappa_x
+ \int_0^{S^\kappa_x}\frac{\kappa^2}{|y|_\kappa+|y|} dy.
\label{2.6a0}
\end{eqnarray}
Thus
\begin{eqnarray}
{  \frac12 (|S^\kappa_x| S^\kappa_x)_x=
 \left(\int_0^{S^\kappa_x}|y|_\kappa dy\right)_x
 - \frac{\kappa^2 S^\kappa_{xx}}{|S^\kappa_x|_\kappa+|S^\kappa_x|}.}
 \label{2.6a1}
\end{eqnarray}
By \eq{lowerBound}  and the Young inequality we obtain from (\ref{5}) and the assumption that $k\le 1$ that
\begin{eqnarray}
 \left|\frac{\kappa^2
 S^\kappa_{xx}}{|S^\kappa_x|_\kappa+|S^\kappa_x|}\right| &\le&
 |\kappa  S^\kappa_{xx}|,\ {\rm thus}\non\\
 \|\kappa S^\kappa_{xx}\|_{L^\frac43(Q_{T_e})} &\le&
 \left(\int_{Q_{T_e}} \left(\kappa^2
 +\kappa|{S^\kappa_{xx}}|^2\right)d xd\tau  \right)^\frac34\le C.
\end{eqnarray}
Combination with \eq{5g}  and  (\ref{2.6a1}) yields
\begin{eqnarray}
 \| (|S^\kappa_x| S^\kappa_x)_x\|_{L^\frac43(Q_{T_e})}
 \le C \left\|\left(\int_0^{S^\kappa_x}|y|_\kappa dy\right)_x\right\|_{L^\frac43(Q_{T_e})}
 + C \|\kappa S^\kappa_{xx}\|_{L^\frac43(Q_{T_e})}  \le C.
 \label{sxsxx}
\end{eqnarray}

It is clear that $ {\overline {|S^\kappa_x| S^\kappa_x} } \le C \int_\Omega |S^\kappa_x| ^2 dx\le C$. Applying
 again the Poincar\'e  inequality to the function $f=|S^\kappa_x| S^\kappa_x$, we arrive at
$$
 \left\|\, |S^\kappa_x| S^\kappa_x\right\|_{L^\frac43(Q_{T_e})} \le C.
$$
Hence this, combined with \eq{sxsxx}, implies that
$$
 \|\, |S^\kappa_x| S^\kappa_x\|_{L^\frac43(0,{T_e};W^{1,\frac43 }(\Omega))} \le C,
$$
one concludes by using the Sobolev embedding theorem that
$$
\|\, |S^\kappa_x| S^\kappa_x\|_{L^\frac43(0,{T_e};L^{\infty}(\Omega))} \le C,
$$
which is
$$
\| S^\kappa_x\|_{L^\frac83(0,{T_e};L^{\infty}(\Omega))} \le C.
$$
This completes the proof of the corollary.

\medskip
To apply some  compactness lemma to the approximate solutions, we need  estimates on the time derivative of the unknown
$S^\kappa$ and also $|S^\kappa_{x}| S^\kappa_{x}$.
\begin{Lemma}\label{L3.3}
 The function $S_t^\kappa$ belongs to $L^\frac43(Q_{T_e})$ and we have the
estimates
\begin{eqnarray}
\label{state0}
\| S^\kappa_t \|_{L^{4/3}(Q_{T_e})} \leq C\,, \\[1ex]
\left\|\left( |S^\kappa_{x}| S^\kappa_{x} \right)_{t} \right\|_{L^1(0,T_e;H^{-2}(\Omega))}
\le C\,.
\label{state1}
\end{eqnarray}
\end{Lemma}

\noindent{\it Proof.} From equation (\ref{m2.2}) and the estimates (\ref{5d}),
 and  (\ref{5}) we immediately see that  $S_t^\kappa\in
L^\frac43(Q_{T_e})$ and that (\ref{state0}) holds. Therefore we only need to
prove the second estimate.

To verify \eq{state1} we must show that there exists a constant $C$,
which is independent of $\kappa$, such that
\begin{eqnarray}
 \left|\left(\left( |S^\kappa_{x}| S^\kappa_{x}\right)_{t},\varphi \right)_{Q_{T_e}}\right|
 \le C\|\varphi\|_{L^\infty(0,T_e;H^2(\Omega))}
 \label{alternative}
\end{eqnarray}
for all functions $\varphi\in L^\infty(0,T_e;H^2_0(\Omega))$. To prove
\eq{alternative}, we first prove that for any $1\ge \delta>0$ there holds
\begin{eqnarray}
 \left|\left(\left(\int_0^{S^\kappa_{x}}
 |y|_\delta dy \right)_{t},\varphi \right)_{Q_{T_e}}\right|
 \le C\|\varphi\|_{L^\infty(0,T_e;H^2(\Omega))}
 \label{alternative1}
\end{eqnarray}
for all functions $\varphi\in L^\infty(0,T_e;H^2_0(\Omega))$. Here
$\delta$ is independent of $\kappa$. Inequality \eq{alternative} is
obtained from this estimate as follows: From $S^\kappa_{x}\in
L^\infty(0,{T_e},L^2(\Omega) ) \subset L^2(Q_{T_e})$,
$S^\kappa_{xt}\in L^2(Q_{T_e})$ and $|\, |y|_\delta - |y|\,|\le
\delta\to 0 $ as $\delta\to 0$ we infer that $\|\,
|S^\kappa_{x}|_\delta - |S^\kappa_{x}|\,\|_{L^\infty(Q_{T_e})} \to 0
$. A straightforward computation yields that
\begin{eqnarray}
 \left(\int_0^{S^\kappa_{x}}
 |y|_\delta dy \right)_{t}
 = |S^\kappa_{x}|_\delta  S^\kappa_{xt}\, .
 \label{formula}
\end{eqnarray}
Therefore, $\left(\int_0^{S^\kappa_{x}} |y|_\delta
dy\right)_t=|S^\kappa_{x}|_\delta S^\kappa_{xt}\to |S^\kappa_{x}| \,
S^\kappa_{xt} $ strongly in $L^2(Q_{T_e})$. Whence, as $\delta\to 0$,
$$
 \left(\left(\int_0^{S^\kappa_{x}} |y|_\delta dy\right)_t,\varphi\right) \to \frac12
 \left( \left( |S^\kappa_{x}| S^\kappa_{x}\right)_{t},\varphi  \right)_{Q_{T_e}}
$$
for all $\varphi\in L^\infty(0,T_e;H^2_0(\Omega))\subset
L^\infty(Q_{T_e}) $. This relation together with \eq{alternative1}
implies \eq{alternative}.

Thus it suffices to prove \eq{alternative1}. To simplify the notations
we define
\begin{eqnarray}
 {\cal R}_\kappa :=  c\nu  |S_{x}^\kappa|_\kappa S^\kappa_{xx} - {\cal F }^\kappa (|S_{x}^\kappa|_\kappa - \kappa).
 \label{Rdef}
\end{eqnarray}
Here ${\cal F }^\kappa = {\cal F }(u^\kappa,u^\kappa_x,S^\kappa, S^\kappa_x)$. Recalling estimate \eq{5max}, we have
\begin{eqnarray}
 |{\cal R}_\kappa |\le C\Big( |S_{x}^\kappa|_\kappa |S^\kappa_{xx}| + (1+| S^\kappa_{x}|)  (|S_{x}^\kappa|_\kappa + \kappa)\Big).
 \label{Rdef1}
\end{eqnarray}
 Multiplying equation (\ref{m2.2}) by $\left(|S^\kappa_{x}|_\delta  \varphi \right)_x, $
integrating the resulting equation with respect to $(t,x)$ over
$Q_{{T_e}}$, using integration by parts for the term with the time
derivative  and noting \eq{formula}, we obtain
\begin{eqnarray}
 0&=&\left(S^\kappa_{t} - {\cal R}_\kappa,\left( |S^\kappa_{x}|_\delta
 \varphi \right)_x\right)_{Q_{T_e}}
 \non\\
 &=& -\left(S^\kappa_{xt}, |S^\kappa_{x}|_\delta
 \varphi \right)_{Q_{T_e}}
 - \left({\cal R}_\kappa, \left( |S^\kappa_{x}|_\delta
 \varphi \right)_x\right)_{Q_{T_e}}\non\\
 &=&-\left(\left(\int_0^{S^\kappa_{x}}
 |y|_\delta  dy \right)_{t},\varphi \right)_{Q_{T_e}}
 -\left({\cal R}_\kappa,
 \left.(|y|_\delta )'\right|_{y=S_x^\kappa} S^\kappa_{xx}\varphi\right)
 -\left({\cal R}_\kappa,  |S^\kappa_{x}|_\delta  \varphi_x\right).\ \
 \label{2.25}
\end{eqnarray}
Remembering that $S^\kappa_{xt}\in L^2(Q_{T_e})$ for any fixed
$\kappa$, we see that the first term in the second equality of
(\ref{2.25}) is properly defined.

To estimate the last two terms on the
right hand side of  inequality \eq{2.25}, we note that there holds
\[
\left| (|y|_\delta )'\right|= \left|\frac{y} {|y|_\delta }\right|\le 1
\quad {\rm and}\quad |y|_\delta \le  |y|+1,
\]
which yields the estimates
\begin{eqnarray}
 && \left|\left({\cal R}_\kappa, \left.(|y|_\delta)'\right|_{y=S_x^\kappa} S^\kappa_{xx}\varphi
 \right)_{Q_{T_e}}\right|
 \le  \left(|{\cal R}_\kappa|, | S^\kappa_{xx}\varphi|
 \right)_{Q_{T_e}}  \non\\
 &\le &  \left(|S^\kappa_{x}|_\kappa |S^\kappa_{xx}|^2,|\varphi| \right)_{Q_{T_e}} +
 \left((1+| S^\kappa_{x}|)  (|S_{x}^\kappa|_\kappa + \kappa),
  |S^\kappa_{xx}\varphi| \right)_{Q_{T_e}}  \non\\
 & \le & C\|\varphi \|_{L^\infty(Q_{T_e})}+ I
 \le C\|\varphi \|_{L^\infty(0,{T_e};H^2(\Omega))}+I,
 \label{2.25z}
\end{eqnarray}
and
\begin{eqnarray}
 \left|\left({\cal R}_\kappa,  |S^\kappa_{x}|_\delta\, \varphi_x\right)_{Q_{T_e}}\right|
 &\le& C \int_{Q_{T_e}}\left(|S^\kappa_{x}| +1\right)|S^\kappa_{x}|_\kappa
  |S^\kappa_{xx}\varphi _x|d(\tau,x)
 \non\\
 & & +\, C \int_{Q_{T_e}}(1+| S^\kappa_{x}|)^2 (|S_{x}^\kappa|_\kappa + \kappa)\,
 \left| \varphi _x\right|d(\tau,x)\non\\
 &=&J_1+J_2.
 \label{2.25y1}
\end{eqnarray}

We estimate $I$ first. Write
\begin{eqnarray}
I &=& C\left(  |S_{x}^\kappa|_\kappa + \kappa   +   | S^\kappa_{x}| \,  |S_{x}^\kappa|_\kappa  +  \kappa | S^\kappa_{x}| ,
  |S^\kappa_{xx}\varphi| \right)_{Q_{T_e}} \non\\
  &=& I_1+I_2+I_3+I_4 .
 \label{Is}
\end{eqnarray}
One needs to estimate $I_i$ ($i=1,2,3,4$). By estimate \eq{5}, $I_1$ can be treated as
\begin{eqnarray}
 I_1 &=& C\left(  |S_{x}^\kappa|_\kappa,   |S^\kappa_{xx}\varphi| \right)_{Q_{T_e}} \non\\
 &\le&  C\int_{Q_{T_e}}  |S^\kappa_{x}|^\frac12_\kappa
 |S^\kappa_{x}|^\frac12_\kappa |S^\kappa_{xx}|\, |\varphi|\, d(t,x)\non\\
 &\le& C\int_0^{{T_e}} \|\,|S^\kappa_{x}|^\frac12_\kappa\|_{L^4(\Omega)}
 \|\,|S^\kappa_{x}|^\frac12_\kappa S^\kappa_{xx}\|\,\| \varphi\|_{L^4(\Omega)}d\tau \non\\
 &\le& C\left(\int_0^{{T_e}}\|\,|S^\kappa_{x}|^\frac12_\kappa S^\kappa_{xx}\|^2d\tau\right)^\frac12
 \left(\int_0^{{T_e}}\| \varphi\|^2_{L^4(\Omega)}d\tau\right)^\frac12\non\\
 &\le& C\|\varphi\|_{L^2(0,T_e;L^4(\Omega) )},
  \label{i1}
\end{eqnarray}
and
\begin{eqnarray}
 I_2 &=& C\left(  \kappa,   |S^\kappa_{xx}\varphi| \right)_{Q_{T_e}} \non\\
 &\le&  C\int_{Q_{T_e}}   \kappa   |S^\kappa_{xx}|\, |\varphi|\, d(t,x)\non\\
 &\le& C  \kappa   \| S^\kappa_{xx}\|_{L^2(Q_{T_e})} \| \varphi\|_{L^2(Q_{T_e})}  \non\\
 &\le& C \|\varphi\|_{L^2(0,T_e;L^2(\Omega) )}.
   \label{i2}
\end{eqnarray}
With the help of \eq{5w} and of the Cauchy-Schwarz inequality, we deal with $I_4$ as follows
\begin{eqnarray}
 I_4 &=& C\left(  \kappa |S_{x}^\kappa| , |S^\kappa_{xx}\varphi| \right)_{Q_{T_e}} \non\\
 &\le&  C\int_0^{T_e} \left( \int_\Omega   |S_{x}^\kappa  S^\kappa_{xx}|^\frac43 dx\right)^\frac34
 \left( \int_\Omega   |\varphi|^4 dx\right) ^\frac14 dt \non\\
 &\le& C  \|\, S_{x}^\kappa S^\kappa_{xx}\|_{L^\frac43(Q_{T_e})} \| \varphi\|_{L^4(Q_{T_e})}  \non\\
 &\le& C \|\varphi\|_{L^4(0,T_e;L^4(\Omega) )}.
   \label{i3}
\end{eqnarray}
The remaining term $I_3$ is the most difficult to evaluate. Making use of estimates \eq{5} and \eq{5f}, we have
\begin{eqnarray}
 & &I_3 = C\left( |S_{x}^\kappa||S_{x}^\kappa|_\kappa , |S^\kappa_{xx}\varphi| \right)_{Q_{T_e}} \non\\
 &\le&  C\| \varphi\|_{L^\infty(Q_{T_e})}
  \int_0^{T_e} \left( \int_\Omega   |S_{x}^\kappa |^2|S_{x}^\kappa|_\kappa dx \right)^\frac12
  \left( \int_\Omega  |S_{x}^\kappa|_\kappa  |S^\kappa_{xx}|^2 dx\right) ^\frac12 dt \non\\
 &\le& C \| \varphi\|_{L^\infty(Q_{T_e})}
  \int_0^{T_e} (\|S_{x}^\kappa\|^\frac12_{L^\infty(\Omega)}+1) \left( \int_\Omega   |S_{x}^\kappa |^2dx \right)^\frac12
  \left( \int_\Omega   |S_{x}^\kappa|_\kappa | S^\kappa_{xx}|^2 dx\right) ^\frac12 dt  \non\\
  &\le& C \| \varphi\|_{L^\infty(Q_{T_e})}
  \left( \int_0^{T_e} (\|S_{x}^\kappa\|^\frac12_{L^\infty(\Omega)}+1)^2 dt \right)^\frac12
  \left( \int_0^{T_e} \int_\Omega   |S_{x}^\kappa|_\kappa | S^\kappa_{xx}|^2 dx dt\right) ^\frac12 dt  \non\\
 &\le& C \| \varphi\|_{L^\infty(Q_{T_e})}.
   \label{i4}
\end{eqnarray}

Next, we consider $J_1, J_2$. The term $J_1$ can be handled as
\begin{eqnarray}
 J_1 &=& C \int_{Q_{T_e}}\left(|S^\kappa_{x}| +1\right)|S^\kappa_{x}|_\kappa
  |S^\kappa_{xx}\varphi _x|d(\tau,x) \non\\
 &\le & C\|\varphi_{x} \|_{L^\infty(Q_{T_e})} \int_0^t \left(\int_\Omega (1+|S^\kappa_{x}| )^4 dx \right)^\frac14
 \left(\int_\Omega (|S^\kappa_{x}|_\kappa |S^\kappa_{xx}|)^\frac43 dx \right)^\frac34 d\tau\non\\
 &\le & C\|\varphi_{x} \|_{L^\infty(Q_{T_e})} \int_0^t \left(1+ \|S^\kappa_{x}\| ^2_{L^\infty(\Omega)}\|S^\kappa_{x}\| ^2 \right)^\frac14
 \left(\int_\Omega (|S^\kappa_{x}|_\kappa |S^\kappa_{xx}|)^\frac43 dx \right)^\frac34 d\tau\non\\
 &\le & C\|\varphi_{x} \|_{L^\infty(Q_{T_e})}
  \left(\int_0^t(1+ \|S^\kappa_{x}\|^2_{L^\infty(\Omega)})d\tau \right)^\frac14
 \left(\int_0^t \int_\Omega (|S^\kappa_{x}|_\kappa |S^\kappa_{xx}|)^\frac43 dx d\tau\right)^\frac34 \non\\
 &\le & C\|\varphi \|_{L^\infty(0,T_e;H^2(\Omega))}.
\end{eqnarray}
Here we used the estimates in \eq{5} and Corollary~\ref{C3.2}, which will also be used to evaluate the term $J_2$. Invoking inequality \eq{Rdef1},
we obtain that
\begin{eqnarray}
 J_2 &\le& C \int_{Q_{T_e}}(1+| S^\kappa_{x}|)^2 (|S_{x}^\kappa|_\kappa + \kappa)\,
 \left| \varphi _x\right|d(\tau,x)\non\\
 &\le& C \|\varphi_{x} \|_{L^\infty(Q_{T_e})} \int_{Q_{T_e}}(1+| S^\kappa_{x}| ^3)  d(\tau,x)\non\\
 &\le& C \|\varphi_{x} \|_{L^\infty(Q_{T_e})} \int_0^{T_e} \left(1+  \| S^\kappa_{x}\|_{L^\infty(\Omega)}
 \int_\Omega  | S^\kappa_{x}| ^2 dx\right) d \tau \non\\
 &\le& C \|\varphi_{x} \|_{L^\infty(Q_{T_e})} \left(1 + \int_0^{T_e}  \| S^\kappa_{x}\|_{L^\infty(\Omega)} d \tau \right)\non\\
 &\le & C\|\varphi\|_{L^\infty(0,T_e;H^2(\Omega))}.
 \label{2.25y}
\end{eqnarray}

Combination of (\ref{2.25}) -- (\ref{2.25y}) and using the Sovolev embedding theorem yield
\begin{eqnarray}
&&\left|\left( \left(\int_0^{S_x^\kappa}|y|_\delta dy
\right)_t,\varphi \right)_{Q_{T_e}}\right|\non\\
&\le& C\left(\|\varphi\|_{L^\infty(0,T_e;H^2_0(\Omega))}+
\| \varphi\|_{L^\infty(Q_{T_e})}+\|\varphi\|_{L^4(0,T_e;L^4(\Omega) )}\right)\non\\[0.2cm]
&\le & C\|\varphi\|_{L^\infty(0,T_e;H^2_0(\Omega))}\,,
\end{eqnarray}
which implies \eq{alternative1} and we complete  the proof.

%
%
\section{Existence of solutions to the phase field model}\label{S4}
We shall use in this section the a priori estimates established in the previous
section to study the convergence of
$(u^\kappa, S^\kappa)$ as $\kappa\to 0$. We shall
show that there is a subsequence, which converges to a weak solution of the
initial-boundary value problem (\ref{Final1}) -- (\ref{Final3}), thereby
proving Theorem~\ref{T1.1}.
\\[1ex]
Note first that the estimates in Corollary~\ref{C3.2}  and Lemma~\ref{L3.3},  the fact that
$\Omega$ is bounded,  and
Poincar\'e's inequality imply
\begin{equation}\label{m4.1}
\|S^\kappa\|_{W^{1,4/3}(Q_{T_e})} \leq C\,,
\end{equation}
for a constant $C$ independent of $\kappa$. Hence, we can select a sequence
$\kappa_n \rightarrow 0$ and a function $S\in W^{1,4/3}(Q_{T_e})$,
such that the sequence $S^{\kappa_n}$, which we again denote by $S^\kappa$,
satisfies
\begin{equation}\label{m4.2}
\| S^\kappa - S \|_{L^{4/3}(Q_{T_e})} \rightarrow 0,\qquad
S^\kappa_x \rightharpoonup S_x\,,\qquad S^\kappa_t \rightharpoonup S_t\,,
\end{equation}
where the weak convergence is in $L^{4/3}(Q_{T_e})$\,.

As usual, since equation (\ref{m2.2}) is nonlinear, the weak convergence
of $S^\kappa_x$ is not enough to prove that the limit function solves this
equation. In the following lemma we therefore show that $S^\kappa_x$ converges
pointwise almost everywhere:
\begin{Lemma}\label{L4.1} There exists a subsequence of $S_x^\kappa$,
we still denote it by $S_x^\kappa$,
such that
\begin{eqnarray}
 \label{m4.3}
 S_x^\kappa\to S_x, & & {  a.e. \ \   in\  \ } Q_{T_e},\\
 \label{m4.4}
 |S_x^\kappa|_\kappa\to |S_x|, & & a.e. \  \   in\  \  Q_{T_e},\\
 \label{m4.5}
 |S_x^\kappa|_\kappa\rightharpoonup
 |S_x|,& &{weakly \ in }\  L^\frac43(Q_{T_e}),\\
 \label{m4.3a}
 \int_0^{S_x^\kappa}  |y| dy \to \frac12 S_x|S_x|, & &  strongly
 \  \  in \  \ {L^\frac43(0,T_e;L^2(\Omega))},\\
 \label{m4.3b}
 \int_0^{S_x^\kappa}  |y|_\kappa dy \to \frac12 S_x|S_x|, & &  strongly
 \  \  in \  \ {L^\frac43(0,T_e;L^2(\Omega))},
\end{eqnarray}
 as $\kappa\to 0$.

\end{Lemma}
The proof is based on  the following two results:
\begin{Theorem}\label{T4.1}
Let $B_0$ be a normed linear space imbedded compactly into
 another normed linear space $B$ which is continuously  imbedded into a
Hausdorff locally convex space $B_1$. Assume that $1\le p<+\infty$,
that $v,v_i \in L^p(0,T_e;B_0)$ for all $i\in \cN$, that the sequence
$\{v_i\}_{i\in \cN}$ converges weakly to $v$ in $L^p(0,T_e;B_0)$ and
that $\{\frac{\partial v_i}{\partial t} \}_{i\in \cN}$ is bounded in
$L^1(0,T_e;B_1)$. Then $v_i$ converges to $v$ strongly in
$L^p(0,T_e;B)$.
\end{Theorem}
\begin{Lemma}\label{L4.2}
Let $(0,T_e)\times \Omega$  be an open
set in $\cR^+\times \cR^n$ and assume that $1<q<\infty$. Suppose that
the functions $g_n, g \in L^q((0,T_e)\times \Omega )$ satisfy
$$
\|g_n\|_{L^q((0,T_e)\times \Omega )}\le C, \ \ g_n\to g \ almost\ everywhere \
in\ (0,T_e)\times \Omega .
$$
Then $g_n$ converges to $g$ weakly in $L^q((0,T_e)\times \Omega )$.
\end{Lemma}
Theorem~\ref{T4.1} is a general version of Aubin-Lions lemma valid under
the weak assumption $\partial_t v_i \in L^1(0,T_e;B_1)$. This
version, which we need here, is proved in Simon \cite{Simon} and in Roub\'ic\v ek
\cite{Roubicek}. A proof of Lemma~{5.2} can be found in
\cite[p.12]{Lions}.
\\[1ex]
{\it Proof of Lemma~\ref{L4.1}:} We choose $p=\frac43$
and
$$
 B_0=W^{1,\frac43}(\Omega),\quad  B=L^2(\Omega),\quad   B_1=H^{-2}(\Omega).
$$
These spaces satisfy the assumptions of the theorem. Since the
estimates (\ref{5d}), (\ref{5g}) and \eq{state1} imply that the sequence
 $\int_0^{S_x^\kappa}  |y| dy$
is uniformly bounded in
$L^p(0,T_e;B_0) $
for $\kappa\to 0$ and $\left(\int_0^{S_x^\kappa}  |y| dy\right)_t $ is
uniformly bounded in $L^1(0,T_e;B_1)$,
it follows from Theorem~3.1 that there is a subsequence, still
denoted by $\int_0^{S_x^\kappa}  |y| dy $, which converges strongly
in $L^p(0,T_e;B) = L^\frac43(0,T_e;L^2(\Omega))$ to a limit function
$G\in L^\frac43(0,T_e;L^2(\Omega)) $.  Consequently, from this
 sequence we can select another
subsequence, denoted in the same way,
which converges almost
everywhere in $Q_{T_e}$. Using that the mapping
$y\mapsto f(y):=\int_0^y|\xi|d\xi =\frac12y|y|$
has a continuous inverse $f^{-1}:\cR\to \cR$, we infer that also the
sequence $S^\kappa_{x} = f^{-1}\left(\int_0^{S_x^\kappa}  |y| dy \right)$
converges pointwisely almost everywhere to $ f^{-1}(G)$ in
$Q_{T_e}$. From the uniqueness of the weak limit we conclude that $
f^{-1}(G)=S_x$ almost everywhere in $Q_{T_e}$.

For the proof of (\ref{m4.3b}) we write
$$
 \int_0^{S_x^\kappa}  |y|_\kappa dy = \int_0^{S_x^\kappa}  |y| dy +
 \int_0^{S_x^\kappa}  (|y|_\kappa - |y|) dy=I_1+I_2.
$$
It is easy to estimate $I_2$ as
$\|I_2\|_{L^2(Q_{T_e})}\le  \|\kappa S_x^\kappa\|_{L^2(Q_{T_e})}
\le C \kappa \|S_{x}^\kappa\|_{L^\infty(0,{T_e};L^2(\Omega))}\le C
\kappa \to 0$.
Therefore, $\int_0^{S_x^\kappa}  |y|_\kappa dy\to \lim_{\kappa\to 0}I_1=\frac12 |S_x |S_x$
strongly in $L^\frac43(0,{T_e};L^2(\Omega)$. This is (\ref{m4.3b}).

To prove (\ref{m4.5}) we note that the estimate
$|S_x^\kappa|_\kappa\le |S_x^\kappa| + \kappa $ and the inequality
(\ref{m4.1}) together imply that the sequence  $|S_x^\kappa|_\kappa
$ is uniformly bounded in $L^\frac43(Q_{T_e})$. Thus, (\ref{m4.5})
is a consequence of (\ref{m4.4}) and Lemma~{5.2}.

\medskip
\noindent{\it Proof of Theorem~1.1:} Define the function $u$ by inserting $S$ into
 \eq{m2.1d} where  $S$ is the limit function of the sequence $S^\kappa$.
We shall prove that $(u, S) $ is a weak solution of problem
(\ref{Final1}) -- (\ref{Final3}).

Remember first that by Lemma~\ref{L3.1} we have $S\in L^\infty(Q_{T_e})$. From this relation, from the above
definition of $u$   we immediately see that $u$  satisfies (\ref{property01}). Observe
next that $\|S^\kappa\|_{L^\infty(0,T_e;H^1_0(\Omega))}\le C,$ by
Lemma~\ref{L3.1} and Sobolev's embedding theorem. This implies $S\in
L^\infty(0,T_e;H^1_0(\Omega))$, since we can select a subsequence of
$S^\kappa$ which converges weakly to $S$ in this space. Thus, $S$
satisfies (\ref{property02}).

Noting that from (\ref{convolution}) and (\ref{m4.2})
\begin{eqnarray}
 \|\chi_\kappa*S^\kappa-S \|_{{L^\frac43}(Q_{T_e})}&\le &
 \|\chi_\kappa*(S^\kappa-S)\|_{{L^\frac43}(Q_{T_e})}+
 \|(S-\chi_\kappa*S) \|_{{L^\frac43}(Q_{T_e})}\non\\
 & \le & \|(S-\chi_\kappa*S) \|_{{L^\frac43}(Q_{T_e})}
 + \|S^\kappa-S \|_{{L^\frac43}(Q_{T_e})}\to 0,
 \label{convergConvolution}
\end{eqnarray}
for $\kappa\to 0$, we conclude easily that the function  $u$  defined in this way satisfy wealkly equation (\ref{Final1}).
 It is thus enough to prove that the equation  (\ref{Final1}) -- (\ref{Final2}) are fulfilled in the weak sense. By definition, these
 equation are   satisfied in the weak sense if the relation
(\ref{definition}) holds. To verify (\ref{definition}) we use that by
construction $ S^\kappa $ solves (\ref{m2.2}). Now we multiply equation (\ref{m2.2}) by a test
function $\varphi \in C_0^\infty((-\infty,T_e)\times\Omega)$ and
integrate the resulting equation over $Q_{T_e}$, then obtain
\begin{eqnarray*}
 0 &=& (S^\kappa_t,\varphi)_{Q_{T_e}} +\left(-c\, \nu |S_x^\kappa |_\kappa  S_{xx}^\kappa
 + {\cal F}^\kappa (|S^\kappa _{x}|_\kappa- \kappa),\varphi \right)_{Q_T}\non\\
 &=& -(S^\kappa_0,\varphi(0))_{\Omega}-(S^\kappa,\varphi_t)_{Q_{T_e}}
 +\left(c\, \nu  \int_0^{S^\kappa_x }|y|_\kappa dy,\varphi_{x}\right)_{Q_{T_e}}\non\\
 & & + \left( {\cal F}^\kappa (|S^\kappa _{x}|_\kappa- \kappa),\varphi \right)_{Q_T}.
\end{eqnarray*}
Equation (\ref{definition}) follows from this relation if we show that
\begin{eqnarray}
 (S^\kappa_0,\varphi(0))_{\Omega} &\to&
 (S_0,\varphi(0))_{\Omega},
 \label{converg0}
\\[1ex]
(S^\kappa,\varphi_t)_{Q_{T_e}} &\to& (S,\varphi_t)_{Q_{T_e}},
 \label{converg1}
\\[1ex]
\left(\int_0^{S^\kappa_x }|y|_\kappa dy,\varphi_x\right)_{Q_{T_e}}
&\to& \left(\frac12|S_x|S_x,\varphi_x\right)_{Q_{T_e}},
 \label{converg2}
\\
\left( {\cal F}^\kappa (|S^\kappa _{x}|_\kappa- \kappa), \varphi\right)_{Q_{T_e}}
&\to& \left( {\cal F} , \varphi\right)_{Q_{T_e}},
 \label{converg3}
\end{eqnarray}
for $\kappa\to 0$. Now, the relation (\ref{converg0}) follows from
(\ref{approximate}), and the relation (\ref{converg1}) is a
consequence of (\ref{m4.2}). By \eq{m4.3b}, one has \eq{converg2}.

To verify (\ref{converg3}) we note that \eq{convergConvolution},
\eq{5f} and the definition of ${\cal F}^\kappa$ yield
\begin{eqnarray}
 \|{\cal F}^\kappa (|S^\kappa _{x}|_\kappa- \kappa)\|_{ L^\frac43(Q_{T_e})}&\le &C,\\[0.2cm]
 {\cal F}^\kappa (|S^\kappa _{x}|_\kappa- \kappa) &\to& {\cal F} |S  _{x}|  ,\ {\rm almost\ everywhere}.
 \label{m2.8a2}
\end{eqnarray}
Then by Lemma~\ref{L4.2},
$$
 {\cal F}^\kappa  (|S^\kappa _{x}|_\kappa- \kappa) \rightharpoonup  {\cal F}  |S^\kappa _{x}| ,
$$
weakly in $L^\frac43(Q_{T_e})$,
which implies (\ref{converg3}). Consequently (\ref{definition}) holds.

It remains to prove that the solution has the regularity properties
 stated  in  (\ref{proper1}) and (\ref{proper2}).
The relation $S_t\in L^\frac43(Q_{T_e})$ is implied by  (\ref{m4.2}).
To verify the second assertion in
(\ref{proper1}), we use estimate (\ref{5f}) to get
$$
\int_0^{T_e}\|S_x^\kappa\|_{L^\infty(\Omega)}^\frac83 dt\le C.
$$
This inequality and $S_x^\kappa\rightharpoonup S_x$ in
$L^\frac83(0,{T_e};L^\infty(\Omega))$ imply $ S_x\in
L^\frac83(0,{T_e};L^\infty(\Omega)).$

To prove (\ref{proper2}), we recall that $\int_0^{S_x^\kappa}|y|_\kappa dy
$ converges to $|S_x|S_x$ strongly in the space
$L^\frac43(0,{T_e};L^2(\Omega))\subset L^\frac43(Q_{{T_e}})$
 and that $\left(\int_0^{S_x^\kappa}|y|_\kappa dy  \right)_x$
 is uniformly bounded in $L^\frac43(Q_{{T_e}})$ for $\kappa\to 0$, by
(\ref{5d}). This together implies that
$\left(|S_x|S_x\right)_x\in L^\frac43(Q_{{T_e}})$.
Finally, to prove the second assertion of  (\ref{proper2}) we choose
a test function $\varphi\in
L^4(0,T_e,W^{1,4}_0(\Omega))$, multiply equation (\ref{m2.2}) by
$-\varphi_x$ and integrate the resulting equation over $Q_{T_e}$ to
  obtain
\begin{equation}
 0=\left(S^\kappa_t-{\cal R}_\kappa,- \varphi_x\right)_{Q_{T_e}}
 =\left(S_{xt}^\kappa,\varphi\right)_{Q_{T_e}}+\left({\cal
    R}_\kappa,\varphi_x
 \right)_{Q_{T_e}},
 \label{last}
\end{equation}
with ${\cal R}_\kappa$ defined in (\ref{Rdef}). Invoking the estimates (\ref{5})  and (\ref{5f}) we deduce
 that
$$
 \|{\cal R}_\kappa\|_{L^\frac43(Q_{T_e})}\le C,
$$
hence equation (\ref{last}) yields
$$
\left(S_{xt}^\kappa,\varphi\right)_{Q_{T_e}}\le
\|{\cal R}_\kappa \|_{L^\frac43(Q_{T_e})}\|\varphi_x\|_{L^4(Q_{T_e})}
\le C\|\varphi\|_{L^4(0,T_e;W^{1,4}_0(\Omega))}\,,
$$
and this means that $S_{xt}^\kappa$ is uniformly bounded in $L^\frac43(0,T_e;
W^{-1,\frac43}(\Omega))$. From this estimate and from
$S_{t}^\kappa\rightharpoonup S_{t}$ in $L^\frac43(Q_{T_e})$ we deduce easily
that $S_{xt}$ belongs to the dual space of $L^4(0,T_e;W^{1,4}_0(\Omega))$,
 which is $L^\frac43(0,T_e;W^{-1,\frac43}(\Omega))$.

\bigskip
\bigskip
\noindent{\bf Acknowledgement.} The authors would like to express their sincere thanks
 to Prof. H.-D. Alber for helpful discussions. The second author of this work has
been partly supported by Grant MTM2008-03541 of
the Ministerio de Educac\'ion y Ciencia of Spain, and by Project
PI2010-04 of the Basque Government.

\end{document}